\documentclass{informspre} 

\OneAndAHalfSpacedXI 

\usepackage{natbib}
\usepackage{tabularx}
\usepackage{algpseudocode}
\usepackage{algorithmicx}
\usepackage[ruled]{algorithm}
\usepackage{caption}
 \usepackage{subcaption}
 \bibpunct[, ]{(}{)}{,}{a}{}{,}%
 \def\bibfont{\small}%

\newcommand{\modif}[1]     {\textcolor[rgb]{0.00,0.00,0.00}{#1}}

\TheoremsNumberedThrough     

\EquationsNumberedThrough    

\begin{document}


\RUNAUTHOR{Zetina, Contreras, and Jayaswal}

\RUNTITLE{Benders Decomposition for Non-Convex Quadratic Facility Location}



\TITLE{An Exact Algorithm for Large-Scale Non-convex Quadratic Facility Location}

\ARTICLEAUTHORS{%
\AUTHOR{Carlos Armando Zetina}
\AFF{IVADO Labs and Canada Excellence Research Chair in Data Science for Real-Time Decision-Making, Polytechnique Montr\'{e}al, Canada, \EMAIL{carlos.zetina@ivadolabs.com} \URL{}}

\AUTHOR{Ivan Contreras}
\AFF{Concordia University and Interuniversity Research Centre on Enterprise Networks, \\Logistics and Transportation (CIRRELT) Montreal, Canada H3G 1M8, \EMAIL{ivan.contreras@cirrelt.ca} \URL{}}

\AUTHOR{Sachin Jayaswal}
\AFF{Indian Institute of Management Ahmedabad, Gujarat, India,  \EMAIL{sachin@iima.ac.in} \URL{}}

} 

\ABSTRACT{
We study a general class of quadratic capacitated $p$-location problems facility location problems with single assignment where a non-separable, non-convex, quadratic term is introduced \modif{in the objective function} to account for the interaction cost between facilities and customer assignments. This problem has many applications in the field of transportation and logistics where its most well-known special case is the single-allocation hub location problem and its many variants. The non-convex, \modif{binary quadratic} program is linearized by applying a reformulation-linearization technique and the resulting continuous auxiliary variables are projected out using Benders decomposition. \modif{The obtained Benders reformulation is then solved using an exact branch-and-cut algorithm that exploits the underlying network flow structure of the decomposed separation subproblems to efficiently generate strong Pareto-optimal Benders cuts.}  Additional enhancements such as a matheuristic, a partial enumeration procedure, and variable elimination tests are also embedded in the proposed algorithmic framework. Extensive computational experiments \modif{on benchmark instances (with up to 500 nodes) and on a new set of instances (with up to 1,000 nodes) of four variants of single-allocation hub location problems confirm the algorithm's ability to scale to large-scale instances.}}
	


\KEYWORDS{quadratic facility location; hub location; \modif{binary quadratic} programs; Benders decomposition;}

\maketitle
%


\section{Introduction}\label{Sec:Intro}\label{S:1}

\textit{Discrete facility location problems} (FLPs) constitute a fundamental class of problems in location science. Given a finite set of customers and a finite set of potential locations for the facilities, FLPs consist of determining the location of facilities to open and the assignment of customers to the selected facilities in order to optimize some cost (or service)-oriented objective. FLPs arise naturally in a wide variety of applications such as supply chain management, emergency and health care systems, humanitarian logistics, telecommunications systems, districting, and urban transportation planning, among others \citep[see,][]{Hamacher2001,Laporte2019}. 

Given the wide range of applications, many variants of FLPs have been studied and they differ with respect to the considered objective, type of demand, network topology, and assignment pattern, among other features. Some of the classical problems involving a cost-oriented objective are \textit{fixed-charge facility location} \citep{Kuehn1963,fernandez15,Fisch17}, \textit{$p$-median} \citep{hakimi,daskin15}, \textit{uncapacitated $p$-location} \citep{Nemhauser1978,wolsey1983,Ortiz2017}, and \textit{single-source capacitated facility location problems} \citep{neebe83,Gadegaard2018}\modif{.} FLPs with service-oriented objectives are \textit{$p$-center} \citep{hakimi,calik15} and \textit{covering location problems} \citep{toregas,Church1974,garcia15}. In such classical FLP variants, the associated objective functions are assumed to be linear functions of the setup cost for opening facilities and/or of the allocation cost of the demand served by the open facilities. 


In this paper, we study a general class of \textit{quadratic capacitated $p$-location problems with single assignments} (QC$p$LPs) which can be defined as follows. Let $N$ be a set of nodes representing both the set of potential facility locations and the set of customers. Let $f_k$ and $b_k$ denote the setup cost and capacity, respectively, of potential facility $k \in N$, and $c_{ik}$ denote the cost of allocating customer $i \in N$ to facility $k \in N$. Let $q_{ikjm}$ represent the interaction cost between the allocation of customer $i \in N$ to facility $k \in N$ and the allocation of customer $j \in N$ to facility $m \in N$. When $i = k$ and $j = m$, $q_{kkmm}$ corresponds to the interaction cost between facilities $k$ and $m$. The QC$p$LP consists of selecting a set of facilities to open, such that no more than $p$ facilities are opened, and of assigning each customer to exactly one open facility while respecting capacity constraints on the amount of demand served by each facility. The objective is to minimize the total setup cost, allocation cost and the (quadratic) interaction cost between facilities and between customer assignments.

The main contributions of this paper are the following. The first is to introduce the QC$p$LPs, which can be stated as a \modif{binary quadratic} program (BQP) having a non-separable, non-convex quadratic function.  QC$p$LPs are quite a general class of problems that can be used to model numerous applications in transportation and telecommunications where a cost (or benefit) captures the interaction between facilities and/or customer assignments. It contains several special cases of hub location problems that are of particular interest such as the \textit{$p$-hub median problem with single assignments} \citep{o1987quadratic,ernst98,Ernst1996}, the \textit{uncapacitated hub location problem with single assignments} \citep{campbell1994integer,Ghaffarinasab2018}, the \textit{capacitated hub location problem with single assignments} \citep{contreras2011b,Meier2018}, and the \textit{star-star hub network design problem} \citep{Helme89,labbe08}. It also has as particular cases some BQPs such as the \textit{quadratic semi-assignment problem} \citep{saito09,Rostami18} and its capacitated variant. In addition, the QC$p$LP can also be used to model situations in which a negative interaction between obnoxious (or undesirable) facilities exist and as a consequence, it is desirable to locate them far from each other. For instance, dangerous plants producing or handling hazardous materials are a source of potential damage to population and their negative impact can be measured as a function of the distances between them \citep{carrizosa15}. A negative interaction among facilities also arises when the facilities compete for the same market and locating facilities as far away from each other is beneficial to mitigate cannibalization and to optimize their competitive market advantage \citep{curtin06,Lei13}. 

The second contribution of this paper is to introduce and computationally compare two strong mixed-integer linear programming (MILP) reformulations of  QC$p$LP. The first is obtained by applying the \textit{reformulation-linearization technique} (RLT) of \cite{adams86} to the original \modif{BQP}. The second reformulation is obtained by projecting out a large number of (continuous) variables from a relaxation of the RLT reformulation using Benders decomposition. Although the resulting Benders reformulation is not as strong as the complete RLT reformulation, its \modif{separation problem has}  a network flow structure, which can be exploited to solve QC$p$LPs efficiently. The third contribution is to present an exact branch-and-cut algorithm based on the Benders reformulation to solve large-scale instances of the QC$p$LP.  We develop several algorithmic refinements to accelerate its convergence. These include: \textit{i)} an exact separation procedure that efficiently generates Pareto-optimal cuts at fractional and integer points by solving transportation problems, \textit{ii)} the use of a matheuristic to generate high quality initial feasible solutions, \textit{iii)} variable elimination tests performed at the root node, and \textit{iv)} a partial enumeration phase to further reduce the size of the problem before branching is performed. To \modif{assess} the performance of our algorithm, we performed extensive computational experiments on several sets of benchmark large-scale instances \modif{of networks} with up to 1,000 nodes for special uncapacitated and capacitated cases previously studied in the literature.


The remainder of the paper is organized as follows. Section \ref{Sec:litrev} provides a succinct review on nonlinear discrete FLPs and positions this paper in view of related work.  Section \ref{Sec:Problem} formally defines the problem and presents a \modif{BQP} formulation. The \modif{RLT} reformulation and the Benders reformulation are presented in Sections \ref{RLTsection}  and \ref{Sec:Benders}, respectively. Section \ref{Sec:Enhance} presents the proposed acceleration techniques and the branch-and-cut algorithm. Section \ref{Sec:CompExp} presents the results of extensive computational experiments performed on a wide variety of \modif{large-scale} instances \modif{of variants of single allocation hub location problems}. Conclusions follow in Section \ref{Sec:Conclusions}.

\section{Literature Review}\label{Sec:litrev}

The study of nonlinear discrete FLPs began with the work of \cite{zangwill68} where production and distribution costs are considered to be separable concave functions of the quantities produced and distributed. In particular, \cite{zangwill68} considers an uncapacitated FLP with arbitrary concave cost functions and presents an exact algorithm that uses a characterization of the extreme points of the feasible region and dynamic programming. \cite{soland74} extends the previous problem to the capacitated case and develops an exact branch-and-bound algorithm in which dual bounds are obtained at the nodes of the enumeration tree by solving a transportation problem. \cite{o1986activity,Kelly86,o1987quadratic} are arguably the first studies dealing with quadratic extensions of FLPs. In particular, \cite{o1987quadratic} introduce a quadratic extension of the $p$-median problem to design hub-and-spoke networks where a set of hub facilities are used as transshipment, consolidation or sorting points to connect a large set of interacting nodes. A non-convex quadratic term arises from the interaction of hub facilities due to the routing of flows between node pairs. The seminal work of \citeauthor{Kelly86} gave rise to a new family of FLPs, denoted as \textit{hub location problems} (HLPs), \modif{that} has since evolved into a rich research area \citep[see,][]{campbell2012twenty,contreras2019}. HLPs lie at the heart of network design \modif{and} planning in transportation and telecommunication systems. Application areas of HLPs in transportation include air freight and passenger travel, postal delivery, express package delivery, trucking, liner shipping, public transportation, and rapid transit systems. Applications of HLPs in telecommunications arise in the design of various digital data service networks. For several particular classes of HLPs, exact solution algorithms have been developed using various integer programming techniques: \textit{i)} branch-and-bound algorithms that use bounding procedures such as combinatorial algorithms \citep{ernst98}, dual-ascent algorithms \citep{canovas07}, and Lagrangean relaxations \citep{contreras11a,Rostami2016,tanash2016,Alibeyg2018}; \textit{ii)} branch-and-price algorithms \citep{contreras2011b,Rothenbcher2016}; \textit{iii)} branch-and-cut algorithms \citep{labbe2005branch,Contreras2012supermodula,Zetina2017TRB,Meier2018}; and \textit{iv)} Benders decomposition algorithms \citep{camargo09,contrerasb,contreras11b,martins15,maheo19}.

\cite{Helme89} study a quadratic extension of the uncapacitated fixed-charge facility location problem arising in the design of satellite communication networks. In this case, demand nodes in the same cluster communicate via terrestrial links whereas different clusters communicate with each other via satellite links using earth stations. A (non-convex) quadratic term is considered in the objective to account for both the within and between cluster traffic costs. \cite{desrochers95} introduce a congested facility location problem in which congestion takes place at facilities and is characterized by convex delay functions that approximate the mean delay experienced by users waiting to be served at a facility. The authors propose a column generation algorithm as a bounding procedure within a branch-and-price algorithm to obtain optimal solutions of the problem. \cite{Fischetti2016} study a congested capacitated facility location problem in which the congestion cost is modeled as a convex quadratic function. The authors present a perspective formulation and an associated Benders reformulation in which all continuous variables are projected out. The resulting convex mixed integer nonlinear program is solved with a branch-and-cut algorithm to obtain optimal solutions for the linear and quadratic capacitated variants. \cite{Harkness03} study a closely related problem in which unit production costs increase once a certain scale of output is reached. Such diseconomies of scale in production are modeled using convex functions to represent situations in which production capacity can be stretched by incurring some additional cost due to overtime, procuring more costly materials, or neglecting equipment maintenance schedules. The authors propose alternative linear integer programming formulations using piecewise linear approximations of the convex functions. \cite{fatma14} focus on a more comprehensive model in which both economies and diseconomies of scale are jointly considered. In the context of production-distribution systems, unit production costs tend to initially decrease as production volume increases, but after some threshold production level, unit costs start to increase due to over-utilization of resources, overtime and facility congestion. Economies and diseconomies of scale are captured by an inverse S-shaped function which is initially concave and then turns convex. The authors present various column generation heuristics to approximately solve this problem. 

\cite{gunluk07} and \cite{gunluk12} study a quadratic extension of the fixed-charge facility location problem in which a separable convex quadratic term is used. In particular, allocation costs are assumed to be proportional to the square of a customer's demand served by an open facility. This problem arises in the design of electric power distribution systems where distribution transformers are located close to customer nodes so as to minimize power loss. Because of cable resistance, power loss in a cable is proportional to the square of the current and the resistance of the cable. \cite{Fisch17} develop a Benders decomposition algorithm based on the perspective formulation of \cite{gunluk12} for solving large-scale instances for both the linear and separable convex quadratic variants. Finally, we refer the reader to \cite{Elhedhli12}, \cite{fatma14}, and \cite{berman2015} for additional references on nonlinear FLPs.

From a methodological perspective, our paper is related to \cite{contrerasb} where an exact algorithm based on Benders decomposition is presented to solve the uncapacitated hub location problem with multiple assignments. Similar to \cite{contrerasb}, we use a Benders reformulation in combination with several algorithmic refinements such as the separation of Pareto-optimal cuts, partial enumeration, and heuristics to obtain optimal solutions to large-scale instances. However, the way these techniques are used in our exact algorithm are significantly different. The Benders reformulation used in \cite{contrerasb} is obtained from a (mixed-integer linear) path-based formulation whereas the one we present in this work is obtained from a relaxation of a strong RLT reformulation of a BQP. \cite{contrerasb} use a standard iterative Benders decomposition algorithm in which at every iteration, a relaxed integer master problem is optimally solved to derive a cut. In our work, we solve the Benders reformulation with a branch-and-cut framework in which Benders cuts are separated at fractional and integer solutions in a single enumeration tree. Moreover, the partial enumeration, heuristic, and Benders cuts of our algorithm exploit the information from fractional solutions whereas in \cite{contrerasb} only integer solutions are used for these techniques. Finally, the separation procedure for generating (non-stabilized) Pareto-optimal cuts presented in \cite{contrerasb} is approximate and can only be used for integer solutions, whereas in our algorithm we use an efficient exact separation routine to generate stabilized Pareto-optimal cuts for both fractional and integer solutions.

\section{Problem Definition}\label{Sec:Problem}


We recall that QC$p$LP seeks to select a set of facilities to open, such that no more than $p$ facilities are opened, and to assign each customer to exactly one open facility while respecting capacity constraints on the amount of demand served by each facility. The objective is to minimize the total setup cost, allocation cost and the (quadratic) interaction cost between facilities and customer assignments. To formulate the problem, for each pair $i,k \in N$, we define binary location/allocation variables $z_{ik}$ equal to one if and only if node $i$ is assigned to facility $k$. When, $i=k$, variable $z_{kk}$ represents whether a facility is installed or not at node $k$. It is also assumed that when $z_{kk}$ is equal to one, node $k$ will be allocated to itself. Hereafter, whenever obvious, the limits on summations and the $\forall$ symbol associated with the constraints are suppressed for ease of notation. The QC$p$LP can be stated as follows \modif{BQP}:
\begin{eqnarray}
    (QF) \quad \min & &
    \sum_{k \in N} f_{k} z_{kk} + \sum_{i \in N}\sum_{k \in N} c_{ik}z_{ik} + \sum_{i \in N}\sum_{i<j}\sum_{k \in N}\sum_{m \in N}q_{ikjm}z_{ik}z_{jm} \label{eq_Q_obj}\\
\mbox{s.t.} &   &  \sum_{k \in N}z_{ik} = 1 \qquad i \in N \label{eq_Q_1}\\
    &  &  z_{ik} \leq z_{kk}  \qquad i,k \in N, i\neq k \label{eq_Q_2}\\
    & & \sum_{k \in N} z_{kk} \leq p \qquad \label{eq_Q_0} \\
	& & \sum_{i\neq k} d_i z_{ik} \leq \bar{b}_{k}z_{kk}  \qquad  k \in N  \label{eq_Q_3} \\
    &  & z_{ik}\in \{0,1\}   \qquad i, k \in N, \label{eq_Q_5}
\end{eqnarray}
where $\bar{b}_{k}=\left(b_{k}-d_k \right)$. The three terms in the objective function \eqref{eq_Q_obj} capture the (linear) setup cost of opening facilities, the (linear) assignment cost of customers to facilities, and the (quadratic) interaction cost  between facilities and customer assignments. Constraint sets \eqref{eq_Q_1} and \eqref{eq_Q_5} ensure that each node is assigned to exactly one facility while constraints \eqref{eq_Q_2} ensure that customers are assigned to open facilities. Constraint \eqref{eq_Q_0} \modif{states} that at most $p$ facilities can be opened. Finally, constraints \eqref{eq_Q_3} are the capacity restrictions on the amount of demand that can be served by each opened facility. Note that these constraints take into account that if facility $k$ is opened, the demand of node $k$ will always be served by such facility. The non-separable, non-convex quadratic term in the objective function in combination with the capacity and single assignment constraints make QC$p$LP a challenging problem to solve. In the following section, we describe the procedure we followed to linearize the quadratic term in order to obtain tight linear MIP reformulations for the QC$p$LP.

\section{RLT-based Reformulations} \label{RLTsection}


A \textit{standard} strategy \citep{Fortet1960, Glover74} to linearizing the bilinear terms $z_{ik}z_{jm}$ of the objective function \eqref{eq_Q_obj} is to replace them by new \textit{continuous} variables $x_{ijkm}$. The nonlinear relation $x_{ikjm} = z_{ik}z_{jm}$ is then imposed via the following sets of linear constraints: 
\begin{align}
& x_{ikjm} \leq z_{ik} &&  i, j, k, m \in N, i<j \label{std1} \\   
& x_{ikjm} \leq z_{jm} &&   i, j, k, m \in N, i<j  \label{std2} \\
& x_{ikjm} \geq z_{ik} + z_{jm} - 1 &&  i, j, k, m \in N, i<j. \label{std3}
\end{align}

Several improved linearization strategies have been proposed for reformulating zero-one quadratic programs as equivalent linear mixed-integer programs, such as \modif{those} described in \citet{Glover75}, \cite{adams86}, \citet{Sherali2007}, \citet{Liberti2007}, and \cite{caprara08}, among others. In particular, the linearization strategy introduced in \cite{adams86} is known to produce tight reformulations and has been generalized and enhanced to design the so-called \modif{reformulation linearization technique RLT} \citep{Sherali1990,sherali13}. The RLT is a powerful technique whose $n$-th hierarchy produces the convex hull of the mixed integer solutions at the expense of having polynomial terms in the resulting reformulation. 

In what follows, we apply the level-1 RLT to QF in order to obtain tight linear MIP reformulations for the QC$p$LP. In the next section, we show how the additional set of continuous variables $x_{ikjm}$ required in a \textit{reduced} reformulation can be projected out using Benders decomposition to obtain an equivalent linear MIP reformulation for the QC$p$LP in the original space of the $z_{ik}$ variables. 



The level-1 RLT applied to \modif{QF} consists of the following steps:
\begin{itemize}
	\item[(i.)] Form $n^3$ constraints by multiplying the $n$ equality constraints \eqref{eq_Q_1} by each $z_{jm}$, $j,m \in N$:
\begin{align*}
&	\sum_{k \in N} z_{ik}z_{jm} = z_{jm}  &&  i, j, m \in N.
\end{align*}

	\item[(ii.)] Form $n^3(n-1)+n^2+n^3$ constraints by multiplying the $n(n-1)+1+n$ inequality constraints \eqref{eq_Q_2}, \eqref{eq_Q_0}, \eqref{eq_Q_3} and the bound constraints by each $z_{jm}$, $j,m \in N$:
\begin{align*}
&   z_{ik}z_{jm} \leq z_{kk}z_{jm}  &&   i,j,k,m \in N, i \neq k \\ 
& \sum_{k \in N} z_{kk}z_{jm} \leq pz_{jm} &&   j,m \in N \\ 
& \sum_{i\neq k} d_i z_{ik}z_{jm} \leq \bar{b}_{k}z_{kk}z_{jm} &&  j,k,m \in N.\\  
& z_{ik}z_{jm} \leq z_{jm} &&  i,j,k,m \in N.\\  
& 0 \leq z_{ik}z_{jm} &&  i,j,k,m \in N.  
\end{align*}

	\item[(iii.)] Form $n^3(n-1)+n^2+n^3$ constraints by multiplying the $n(n-1)+1+n$ inequality constraints \eqref{eq_Q_2}, \eqref{eq_Q_0}, and \eqref{eq_Q_3} by each $\left( 1-z_{jm}\right)$, $j,m \in N$. Note the bound constraints are ommited as they give meaningless inequalities.:
	\begin{align*}
	&   z_{ik}\left( 1-z_{jm}\right) \leq z_{kk}\left( 1-z_{jm}\right)  &&   i,j,k,m \in N, i \neq k \\ 
	& \sum_{k \in N} z_{kk}\left( 1-z_{jm}\right) \leq p\left( 1-z_{jm}\right) &&   j,m \in N \\ 
	& \sum_{\substack{i \in N \\ i\neq k}} d_i z_{ik}\left( 1-z_{jm}\right) \leq \bar{b}_{k}z_{kk}\left( 1-z_{jm}\right) &&  j,k,m \in N. 
	\end{align*}

	\item[(iv.)] Linearize the above sets of constraints and objective \eqref{eq_Q_obj} by substituting $x_{ijkm}=z_{ik}z_{jm}$, for each $i,j,k,m \in N, i<j$, and add $x_{ikjm} \geq 0$ to obtain the following linear MIP reformulation:
\end{itemize}
\begin{eqnarray}
 (RL_1) \quad  \min & &
    \sum_{k \in N}f_{k} z_{kk} +\sum_{i \in N}\sum_{k \in N} c_{ik}z_{ik} + \sum_{i \in N}\sum_{i<j}\sum_{k \in N}\sum_{m \in N}q_{ikjm}x_{ikjm} \notag\\
\mbox{s.t.} &  & \eqref{eq_Q_1}-\eqref{eq_Q_5} \notag \\
& &   \sum_{k \in N} x_{ikjm} = z_{jm} \qquad   i,j,m \in N, i<j \label{eq_SK_1}\\
&    &   \sum_{k \in N} x_{jmik} = z_{jm}  \qquad  i, j, m\in N, i>j\label{eq_SK_2}\\
&    &   x_{ikjm} \leq x_{kkjm}  \qquad  i,j,k,m \in N, i \neq k, i<j, k<j \label{eqRLT2a}\\
&    &   x_{ikjm} \leq x_{jmkk}  \qquad   i,j,k,m \in N, i \neq k, i<j, k>j \label{eqRLT2b}\\
&    &   x_{jmik} \leq x_{kkjm} \qquad    i,j,k,m \in N, i \neq k, i>j, k<j \label{eqRLT2c}\\
&    &   x_{jmik} \leq x_{jmkk} \qquad    i,j,k,m \in N, i \neq k, i>j, k>j \label{eqRLT2d}\\
&    &   x_{ikjm} \leq z_{jm}  \qquad  i,j,k,m \in N, i \neq k, i<j, k<j \label{eqRLTbound}\\
    & & \sum_{k < j} x_{kkjm} + \sum_{k > j} x_{jmkk} \leq pz_{jm} \qquad  j,m \in N \label{eqRLT1} \\
   & & \sum_{i<j} d_i x_{ikjm} + \sum_{i>j} d_i x_{jmik} \leq \bar{b}_{k}x_{kkjm} \qquad  j,k,m \in N, k < j  \label{eqRLT3a} \\
   & & \sum_{i<j} d_i x_{ikjm} + \sum_{i>j} d_i x_{jmik} \leq \bar{b}_{k}x_{jmkk} \qquad  j,k,m \in N, k > j  \label{eqRLT3b} \\
   & & z_{ik} - x_{ikjm} \leq z_{kk} - x_{kkjm}  \qquad   i,j,k,m \in N, i \neq k, i<j, k<j \label{eqRLT5a}\\
   & & z_{ik} - x_{ikjm} \leq z_{kk} - x_{jmkk}  \qquad  i,j,k,m \in N, i \neq k, i<j, k>j \label{eqRLT5b}\\
   & & z_{ik} - x_{jmik} \leq z_{kk} - x_{kkjm}  \qquad   i,j,k,m \in N, i \neq k, i>j, k<j \label{eqRLT5c}\\
   & & z_{ik} - x_{jmik} \leq z_{kk} - x_{jmkk}  \qquad   i,j,k,m \in N, i \neq k, i>j, k>j \label{eqRLT5d}\\
   & & \sum_{k \in N} z_{kk} - \sum_{k < j} x_{kkjm} - \sum_{k > j} x_{jmkk} \leq p\left( 1-z_{jm}\right) \qquad   j,m \in N \label{eqRLT4} \\
   & & \sum_{i \neq k} d_i z_{ik} - \sum_{i < j} d_ix_{ikjm} - \sum_{i > j}d_i x_{jmik} \leq \bar{b}_{k}\left(z_{kk}-x_{kkjm}\right) \ \ j,k,m \in N, k < j \label{eqRLT6a} \\
   & & \sum_{i \neq k} d_i z_{ik} - \sum_{i < j} d_ix_{ikjm} - \sum_{i > j} d_ix_{jmik} \leq \bar{b}_{k}\left(z_{kk}-x_{jmkk}\right) \ \  j,k,m \in N, k > j \label{eqRLT6b}  \\
& & x_{ikjm} \geq 0 \qquad  i, j, k, m \in N, i<j. \label{eq_SK_3}
\end{eqnarray}

The RLT reformulation of \cite{adams86} contains the standard linearization constraints \eqref{std1}-\eqref{std3}, which are obtained after multiplying $\left( 1-z_{jm}\right)$, $j,m \in N$, by each bound constraint $0\leq z_{ik} \leq 1$, $i,k \in N$. However, these constraints are \modif{implied} by constraints \eqref{eq_SK_1}-\eqref{eq_SK_2} and thus, there is no need to add them to $RL_1$ \modif{\citep{adams86}}. Moreover, an important consequence of this implication is that constraints \eqref{eq_SK_1}-\eqref{eq_SK_2} are sufficient to provide a valid linear MIP formulation of QC$p$LP. That is, constraints \eqref{eqRLT2a}-\eqref{eqRLT6b}, although useful for strengthening the linear programming (LP) relaxation, are actually redundant for the description of the set of feasible solutions to QC$p$LP in the extended space of the $(z,x)$ variables. Therefore, the QC$p$LP can be stated as the following \textit{reduced} linear MIP:
\begin{eqnarray}
(RL_2) \quad  \min & &
\sum_{k \in N}f_{k} z_{kk} +\sum_{i \in N}\sum_{k \in N} c_{ik}z_{ik} + \sum_{i \in N}\sum_{i<j}\sum_{k \in N}\sum_{m \in N}q_{ikjm}x_{ikjm} \notag\\
\mbox{s.t.} &  & \eqref{eq_Q_1}-\eqref{eq_Q_5}, \eqref{eq_SK_1}, \eqref{eq_SK_2}, \eqref{eq_SK_3}. \notag 
\end{eqnarray}

Although both reformulations $RL_1$ and $RL_2$ contain $n^2$ binary variables and $n^3(n-1)/2$ continuous variables, the number of additional constraints required in $RL_1$ is about $n^3(n-1) + 3n^2(n-1)/2 + 2n^2$, whereas in $RL_2$ is only $n^2(n-1)/2$.

We next provide the results of computational experiments to compare the quality of the LP relaxation bounds obtained with $RL_1$ and $RL_2$, as well as with six additional relaxations, denoted as $RL_3, \dots, RL_8$, to \modif{assess} the marginal contribution to the improvement of the LP bounds when adding specific subsets of constraints of $RL_1$ independently. Table \ref{table1} provides the configurations of the considered relaxations of the RLT-based reformulation.

\begin{table}[htbp]
	\TABLE
	{Configuration of considered relaxations of the RLT-based reformulation.\label{table1}}
	{\begin{tabular}{|c|l|}
			\hline \up
			Formulation & Combination of RLT inequalities \\	
			\hline \up
			$RL_1$ & Full RLT reformulation\\
			$RL_2$ & Only RLT for assignment constraints: (10)--(11)\\
			$RL_3$ & $RL_2$ + RLT for linking constraints: (12)--(15)  \\
			$RL_4$ & $RL_2$ + RLT for linking constraints: (19)--(22)  \\
			$RL_5$ & $RL_2$ + RLT for all linking constraints: (12)--(15), (19)--(22)  \\
			$RL_6$ & $RL_2$ + RLT for capacity constraints: (17)--(18)  \\
			$RL_7$ & $RL_2$ + RLT for capacity constraints: (24)--(25)  \\
			$RL_8$ & $RL_2$ + RLT for all capacity constraints: (17)--(18), (24)--(25)   \\
			\hline
	\end{tabular}}
	{}
\end{table}%

Tables \ref{table2} and \ref{table3} provide a comparison of the obtained \%LP gap and computational time (in seconds) needed to solve the LP relaxation, respectively, of each of the considered relaxations. For this experiment, we focus on one well-known particular case of the QC$p$LP, the capacitated hub location problem with single assignments (CHLPSA), and use a set of small-size instances of the AP data set considering $n \in \left\lbrace 20, 25, 40 \right\rbrace$ (see Section \ref{Sec:CompExp} for a detailed description of the instances and the computational setting used throughout our experiments). The entry N.A. corresponds to the cases in which CPLEX failed to solve the associated LP relaxation due to numerical errors encountered during the solution process.

\begin{table}[htbp]
	\TABLE
	{Comparison of \%LP gaps for different RLT-based reformulations for the CHLPSA. \label{table2}}
	{\begin{tabular}{|c|r|r|r|r|r|r|r|r|}
			\hline \up
			& \multicolumn{8}{c|}{\% LP gap} \\			
			\hline \up
			Instance		&	$RL_2$		&	$RL_3$		&	$RL_4$		&	$RL_5$		&	$RL_6$		&	$RL_7$		&	$RL_8$		&$RL_1$	\\
			\hline \up
			20LT 	&	1.59	&	1.30	&	1.30	&	1.30	&	1.49	&	1.00	&	0.97	&	0.37	\\
			20TT 	&	1.64	&	1.63	&	1.63	&	1.63	&	1.61	&	1.06	&	1.03	&	0.90	\\
			25LL 	&	0.17	&	0.17	&	0.17	&	0.17	&	0.17	&	0.15	&	0.15	&	0.10	\\
			25LT 	&	1.88	&	1.88	&	1.88	&	1.88	&	1.88	&	1.21	&	1.18	&	0.62	\\
			25TT 	&	3.34	&	3.00	&	3.00	&	3.00	&	3.27	&	3.01	&	N.A.		&	N.A.	\\
			40LT 	&	0.76	&	0.76	&	0.76	&	0.76	&	0.75	&	0.16	&	0.14	&	N.A.	\\
			40TL 	&	0.03	&	0.03	&	0.03	&	0.03	&	0.03	&	N.A.		&	N.A.		&	0.00	\\
			\hline
	\end{tabular}}
	{}
\end{table}%

\begin{table}[htbp]
	\TABLE
	{Comparison of computational time needed to solve the LP relaxations of different RLT-based reformulations for the CHLPSA. \label{table3}}
	{\begin{tabular}{|c|r|r|r|r|r|r|r|r|}
			\hline \up
			& \multicolumn{8}{c|}{CPU time (seconds)} \\			
			\hline \up
			Instance		&	$RL_2$		&	$RL_3$		&	$RL_4$		&	$RL_5$		&	$RL_6$		&	$RL_7$		&	$RL_8$		&	$RL_1$	\\
			\hline \up
			20LT	&	13	&	172	&	208	&	191	&	158	&	391	&	772	&	39,326	\\
			20TT 	&	11	&	185	&	207	&	217	&	146	&	442	&	834	&	60,384	\\
			25LL 	&	154	&	2,403	&	3,652	&	3,766	&	1,392	&	5,992	&	9,234	&	51,527	\\
			25LT 	&	76	&	923	&	1,338	&	1,339	&	1,059	&	1,904	&	7,817	&	4,736	\\
			25TT 	&	74	&	863	&	1,395	&	1,434	&	1,236	&	3,456	&	N.A.	&	N.A.	\\
			40LT 	&	8,152	&	112,210	&	191,358	&	177,974	&	186,659	&	814,479	&	782,425	&	N.A.	\\
			40TL 	&	9,184	&	141,839	&	239,990	&	233,689	&	183,850	&	N.A.	&	N.A.	&	19,410	\\
			\hline
	\end{tabular}}
	{}
\end{table}%

From Tables \ref{table2} and \ref{table3}, we note that although there is a clear improvement in the quality of the LP bound when solving $RL_1$, the substantial increase in the CPU time required to solve the associated LP when compared to $RL_2$, does not justify its use. The same is true for the rest of the relaxations ($RL_3$ to $RL_8$). Moreover, adding constraints \eqref{eqRLT2a}-\eqref{eqRLT6b} breaks the decomposability and the network flow structure of the separation problem associated with the Benders reformulation of $RL_2$. For these reasons, we next present a Benders reformulation only for $RL_2$. 




\section{Benders Reformulation}\label{Sec:Benders}

Benders decomposition is a well-known partitioning method applicable to MILPs \citep{Benders62}. In particular, it reformulates an MILP by projecting out a set of complicating \modif{continuous} variables to obtain a formulation with fewer variables but typically with a huge  number of constraints, which can be separated efficiently via the solution to an LP subproblem, known as the \textit{dual subproblem} (DSP). These new constraints are usually referred to as Benders cuts and involve only the variables kept in the reduced problem, plus one additional continuous variable. Given that only a small subset of these constraints are usually active in an optimal solution, a natural relaxation is obtained by dropping most of them and generating them on the fly as needed within a cutting plane algorithm. 

Let $Z$ denote the set of vectors $z$ satisfying constraints \modif{\eqref{eq_Q_1}}--\eqref{eq_Q_5}. For any fixed $\bar{z} \in Z$, the \textit{primal subproblem} (PS) in the space of the $x$ variables is
\begin{align}
(PS) \quad \min \quad & \sum_{i \in N}\sum_{j>i}\sum_{k \in N}\sum_{m \in N} q_{ikjm} x_{ikjm} \notag\\
\mbox{s.t.} \quad
&   \sum_{k \in N} x_{ikjm} = \bar{z}_{jm} &&   i,j,m \in N, i<j \label{eq_PS_1}\\
&   \sum_{m \in N} x_{ikjm} = \bar{z}_{ik} &&   i,j,k \in N, i<j \label{eq_PS_2}\\
& x_{ikjm} \geq 0 && i,j,k,m \in N, i<j. \label{eq_PS_3}
\end{align}

Note that the indices $i,j$ and $k,m$ of constraints \eqref{eq_SK_2} have been swapped in \eqref{eq_PS_2} to better highlight the decomposability of the problem, where subproblems consist of several independent transportation problems. In particular, PS can be decomposed into $n(n-1)/2$ problems PS$_{ij}$, one for each pair $(i,j) \in N \times N, i<j$. Therefore, we can construct the corresponding dual subproblem (DS$_{ij}$) for each $(i,j)$.  Let $\alpha$ and $\beta$ denote the dual variables of constraints \eqref{eq_PS_1} and \eqref{eq_PS_2}, respectively. For each $(i,j) \in N \times N, i<j$, the corresponding DS$_{ij}$ can be stated as follows:
\begin{align}
(DS_{ij}) \quad \text{maximize} \quad & \sum_{m \in N}\bar{z}_{jm} \alpha_{ijm} + \sum_{k \in N}\bar{z}_{ik} \beta_{ijk}  \label{eq_DSij_obj}\\
\mbox{s.t.} \quad
&   \alpha_{ijm} + \beta_{ijk} \leq  q_{ikjm} &  k,m  \in N. \label{eq_DSij_1}
\end{align}

Given the network flow structure of each PS$_{ij}$, a sufficient condition for feasibility to PS$_{ij}$ is $\sum_{k\in N}\bar{z}_{ik}=\sum_{m \in N}\bar{z}_{jm}$, for each $(i,j)\in N \times N, i<j$. This is actually guaranteed by constraints \eqref{eq_Q_1} of $Z$ and as a consequence, PS$_{ij}$ is always feasible. Therefore, DS$_{ij}$ has always at least one optimal solution with a finite objective value. In turn, this implies there is no need to incorporate feasibility cuts to the resulting Benders reformulation. For $i,j \in N, i<j$, let $EP_{ij}$ denote the set of extreme points of the polyhedron of DS$_{ij}$. The optimal value of each  DS$_{ij}$ is then equal to 

$$\max_{(\alpha, \beta) \in EP_{ij}} \sum_{m \in N}\bar{z}_{jm} \alpha_{ijm} + \sum_{k \in N}\bar{z}_{ik} \beta_{ijk}. $$

Introducing the extra continuous variable $\eta$ for the overall interaction cost, the \textit{Benders reformulation} (BR) associated with $RL_2$ is
\begin{align}
(BR) \quad \min\quad & \sum_{k \in N}f_{k} z_{kk} + 
\sum_{i \in N}\sum_{k \in N}c_{ik} z_{ik} + \eta  \label{eq_BR_obj_agg}\\  
\text{subject to} \quad & \eqref{eq_Q_1}-\eqref{eq_Q_5} \notag\\
\quad  & \eta \geq \sum_{i \in N} \sum_{k \in N} \left(\sum_{j<i} \alpha_{jik} + \sum_{j>i} \beta_{ijk}\right) {z}_{ik}  && (\alpha, \beta) \in EP, \label{eq_BR_cuts_agg}
\end{align}
where $EP$ is the Cartesian product of the sets of extreme points $EP_{ij}$, $i,j \in N, i<j$ and constraints \eqref{eq_BR_cuts_agg} are the so-called Benders \textit{optimality} cuts. We note that BR contains only the original variables $z$ and one additional continuous variable. 

Whereas it has been empirically shown \citep{Magn81} that the disaggregated form of \eqref{eq_BR_cuts_agg}, \modif{with one for each subproblem}, may lead to an improved performance of the Benders decomposition algorithm for some classes of problems, our preliminary computational experiments showed that not to be the case for our problem, since the need to add $n(n-1)/2$ constraints at every iteration leads to a problem that increases too rapidly in size, especially for medium to large-size instances with $n>200$. As shown in \modif{the} next section, the proposed algorithm does however exploit the decomposability of the dual subproblem leading to a reduction of the required CPU time to prove optimality. 
\section{An Exact Algorithm for QC$p$LP}\label{Sec:Enhance}


In this section, we present an exact branch-and-cut algorithm based on BR to solve QC$p$LPs. The \textit{standard} Benders decomposition algorithm is an iterative procedure in which at every iteration a relaxed integer master problem, containing only a small subset of Benders cuts, is optimally solved to obtain a dual bound. The dual subproblem is then solved to obtain a primal bound and to determine whether additional Benders cuts are needed in the relaxed master problem. If needed, these cuts are added to the master problem and solved again. This iterative procedure is repeated until the convergence of the bounds is attained, if an optimal solution exists. One of the major drawbacks of this \textit{standard} approach is the need to solve an integer master problem at each iteration. To overcome this difficulty, \modif{\textit{recent}} implementations of Benders decomposition have considered the solution of the Benders reformulation with a standard branch-and-cut framework, in which Benders cuts are separated not only at nodes with integer solutions but also at the nodes with fractional solutions of a single enumeration tree \citep[see, for instance][]{Fisch17,Ortiz2017,Zetina2018}. We use this approach to develop an exact algorithm for QC$p$LP. 

In addition, we use the following strategies to speed up the convergence of our branch-and-cut algorithm: i) we exploit the structure of the subproblem to generate Pareto-optimal cuts by efficiently solving network flow problems, ii) we use a simple but effective stabilization procedure for generating cuts, iii) we use a matheuristic to generate high quality solutions, iv) we apply variable elimination tests at the root node, and v) we perform a partial enumeration phase to permanently fix location variables $z_{kk}$ to either one or zero before exploring the enumeration tree.

\subsection{Stabilized Pareto-optimal Cuts}\label{Subsec:Pareto}


It is well known that the selection of Benders cuts plays an important role in the overall convergence of Benders decomposition algorithms. \cite{Magn81} proposed a procedure for obtaining Pareto-optimal cuts, that is, cuts that are not dominated by any other cut. In particular, \cite{Magn81} define the notion of cut dominance as follows. Given two cuts defined by dual solutions $u$ and $u'$ of the form $\theta\geq f(u) +zg(u)$ and $\theta\geq f(u') +zg(u')$, respectively, the cut defined by $u$ dominates the cut defined by $u'$ if and only if $f(u) +zg(u)\geq f(u') +zg(u')$ with strict inequality holding for at least one feasible point $z$ of MP. If a cut defined by $u$ is not dominated by any other cut, then it is a Pareto-optimal Benders cut.

To obtain Pareto-optimal cuts, \cite{Magn81} propose solving an additional linear program similar to DS. It is parameterized by a \textit{core point} of the set $Z$, which is a point in the relative interior of its convex hull. Let $z^0$, $\bar{z}$, and $\Gamma(\bar{z})$ denote a given core point, the current (possibly fractional) solution, and the optimal solution value of DS$_{ij}$, respectively. To determine whether there exists a Pareto-optimal cut or not that is violated by the point $\bar{z}$, for each $i,j \in N, i<j$, we solve the following Pareto-optimal subproblem:
\begin{align}
(DPO_{ij})\quad \text{maximize}\quad &
\sum_{m \in N}z^0_{jm} \alpha_{ijm} + \sum_{k \in N}z^0_{ik} \beta_{ijk}  \label{eq_PODSP_obj}\\
\mbox{s.t.} \quad 
&\sum_{m \in N}\bar{z}_{jm} \alpha_{ijm} + \sum_{k \in N}\bar{z}_{ik} \beta_{ijk} \geq \Gamma(\bar{z}) \label{eq_PODSP_0} &&\\
& \alpha_{ijk} + \beta_{ijm} \leq  q_{ikjm} & & k,m  \in N, \label{eq_PODSP_1}
\end{align}
where constraints \eqref{eq_PODSP_0} guarantee that the optimal solution to $DPO_{ij}$ is contained in the set of optimal solutions to the original $DP_{ij}$.


Let $\delta$ be the dual variable associated with constraints \eqref{eq_PODSP_0} and $x_{ikjm}$ be those of constraints \eqref{eq_PODSP_1}. Dualizing DPO$_{ij}$ we obtain the following linear program:
\begin{align}
(PPO_{ij})\quad \min\quad &
\sum_{k \in N}\sum_{m \in N} q_{ikjm}x_{ikjm}-\Gamma(\bar{z}) \delta \label{eq_POPSP_obj}\\
\mbox{s.t.} \quad 
& \sum_{k \in N}x_{ikjm}=z_{jm}^0+\bar{z}_{jm}\delta& & m\in N \label{eq_POPSP_1}\\
& \sum_{m \in N}x_{ikjm}=z_{ik}^0+\bar{z}_{ik}\delta& & k\in N \label{eq_POPSP_2}\\
& x_{ikjm}\geq 0& & k, m\in N \\
& \delta^{ij}\geq 0. & &
\end{align}

Given that $\delta$ affects the right-and-side of flow constraints \eqref{eq_POPSP_1} and \eqref{eq_POPSP_2}, this problem corresponds to a \textit{parametric transportation problem}.  Given that $\sum_{m \in N} \bar{z}_{jm} = \sum_{k \in N} \bar{z}_{ik} = 1$,  PPO$_{ij}$ can be interpreted as a problem where a rebate of $\Gamma$ is given for each unit of additional supply/demand shipped. Similar to the Pareto-optimal problem presented in \cite{Magn86} for uncapacitated multi-commodity network design problems, rather than performing a parametric analysis on PPO$_{ij}$ to determine its optimal solution, we use the following result.

\begin{proposition}
Any value of $\delta \geq 1$ is optimal to PPO$_{ij}$. 
\end{proposition}
\proof{Proof}
The total demand $\sum_{k \in N}\left( z_{ik}^0+\bar{z}_{ik}\delta\right)$ must flow via some subset of arcs $(k,m) \in N\times N$. At most $\sum_{k \in N}z_{ik}^0 = 1$ of this flow can go through demand nodes $j \in N$ having $\bar{z}_{ik}=0$. Any additional flow must use arcs incident to demand nodes $j \in N$ such that $\bar{z}_{ik}>0$, and thus, the marginal flow cost will be precisely $\Gamma(\bar{z})$. Therefore, any value of $\delta \geq 1$ must be optimal for PPO$_{ij}$.
\endproof

An important consequence of this result is that there is no longer a  need to know the value $\Gamma(\bar{z})$ to select $\delta$. Therefore, the complexity for generating a Pareto-optimal cut is the same as generating a standard Benders optimality cut, i.e, solving $n(n-1)/2$ transportation problems, which can be efficiently solved by using the network simplex algorithm.


In our algorithm, we consider the following family of core points of $Z$ for $p \geq 2$:
%

$$z^0_{ik}(\tilde{H}) = 
\begin{cases}
\frac{p}{|\tilde{H}|}-\epsilon \qquad \quad \ \ i\in \tilde{H}, k \in  \tilde{H}, i = k, \\
\frac{1-\left(\frac{p}{|\tilde{H}|} -\epsilon\right)}{|\tilde{H}|-1} \qquad i\in \tilde{H}, k \in  \tilde{H}, i \neq k, \\
\frac{1}{|\tilde{H}|} \qquad \qquad \quad i\in N\backslash\tilde{H}, k \in  \tilde{H}, i \neq k, \\
\end{cases}
$$
where $0 < \epsilon < \min \left\lbrace \frac{1}{|\tilde{H}|}, \frac{p-1}{|\tilde{H}| - 2} \right\rbrace $, and $\tilde{H} \subseteq N$ denote the current set of candidate facilities to open at a given node of the enumeration tree. 

Although the solution to the Pareto-optimal subproblem $PPO_{ij}$ guarantees that the obtained cut will be non-dominated, in practice some non-dominated cuts may be more useful than others. That is, a core \modif{point} selection strategy is needed to further improve the convergence of the cutting plane algorithm used to solve the LPs at the nodes of the enumeration tree. 

\cite{Ortiz2017}, \cite{Zetina17}, and \cite{Zetina2018} present and computationally compare several static and dynamic core point selection strategies to separate Pareto-optimal cuts for multi-level facility location and multi-commodity network design problems. An interesting observation from these studies is that the strategy providing on average the best results is different in each of these works. These strategies can actually be seen as stabilization procedures for generating effective non-dominated cuts. Stabilization is frequently used to improve the convergence of column generation and cutting plane algorithms needed to solve Dantzig-Wolfe decompositions and Lagrangean relaxation \citep{benameur07, Fisch17, pessoa18}.

After performing extensive preliminary experiments, we note that an efficient core point selection strategy is to dynamically update the separation point $\hat{z}^{'}$ at each iteration by considering a convex combination of the separation point $\hat{z}$ considered in the previous iteration and the current point $\bar{z}$ as follows:
\begin{equation}
\label{updatecore}
\hat{z}^{'} = \phi \hat{z} + (1-\phi)\bar{z},
\end{equation}
where $\phi\in \mathbb{R}$, and $0 < \phi < 1$. We initialize the algorithm by using $\hat{z} = z^0(N)$ at the first iteration of the cutting plane algorithm. Each time the reduction tests of Section \ref{Subsub:Elim} modify the current set of candidate facilities to open  $\tilde{H} \subseteq N$, we reinitialize $\hat{z}$ as $z^0(\tilde{H})$. 



\subsection{A Matheuristic for QC$p$LP}\label{Subsec:Heuristic}
When using branch-and-cut algorithms, it is important to find high quality feasible solutions early on in the process. This leads to smaller search trees since they provide better bounds for pruning and a guide for selecting variables to branch on. Finding near-optimal solutions in a preprocessing stage can be used to perform variable elimination tests that reduce the size of the formulation to be solved \citep{contrerasb,contreras2011b,Alibeyg2018}. 

In this section, we present a matheuristic that exploits the information generated during the solution of the root node to effectively explore the solution space. It consists of two phases. The first is a \textit{constructive heuristic} in which the support of the LP relaxation of the Benders reformulation is used to build a reduced linear integer relaxation of $QF$, which is then solved with a general purpose solver to construct an initial feasible solution. The second is a \textit{local search heuristic} in which several neighborhoods are systematically explored to improve the initial solution obtained during the constructive phase. One of such neighborhoods is efficiently explored by solving an MILP with CPLEX. We next describe in detail each of these phases. 



\subsubsection{Constructive Heuristic}\label{Subsub:MILPRelax}

Let $\bar{z}^t$ denote the solution at iteration $t$ of the LP relaxation of the Benders reformulation at the root node of the enumeration tree and let $H(\bar{z}^t) = \left\lbrace k \in N: \bar{z}^t_{kk} > 0 \right\rbrace $ correspond to the set of (partially) opened hub facilities in  $\bar{z}^t$. We can obtain a feasible solution to QC$p$LP by solving the following reduced (linear) relaxation of SQFLP: 
\begin{eqnarray}
(RLF) \quad \min\quad & &
\sum_{k \in H(\bar{z}^t)} f_{k} z_{kk} + 
\sum_{i \in N}\sum_{k \in H(\bar{z}^t)} c_{ik}z_{ik} \label{eq_L_obj}\\
\mbox{s.t.} \quad & &    \sum_{k \in H(\bar{z}^t)}z_{ik} = 1 \qquad i \in N \label{eq_Q_1red}\\
&  &  z_{ik} \leq z_{kk}  \qquad i \in N, k \in H(\bar{z}^t), i\neq k \label{eq_Q_2red}\\
& & \sum_{k \in H(\bar{z}^t)} z_{kk} \leq p \qquad \label{eq_Q_0red} \\
& & \sum_{i\neq k} d_i z_{ik} \leq \bar{b}_{k}z_{kk}  \qquad  k \in H(\bar{z}^t)  \label{eq_Q_3red} \\
&  & z_{ik}\in \{0,1\}   \qquad i \in N, k \in H(\bar{z}^t). \label{eq_Q_5red}
\end{eqnarray}

Note that the quadratic term of objective function \eqref{eq_Q_obj} of QF has been relaxed. Whenever $|H(\bar{z}^t)| \ll |N|$, RLF will be substantially smaller than QF. Given that any feasible solution to \eqref{eq_Q_1red}-\eqref{eq_Q_5red} is also feasible for \eqref{eq_Q_1}-\eqref{eq_Q_5}, we can compute a valid upper bound on the optimal solution of QC$p$LP by simply evaluating objective \eqref{eq_Q_obj} using the optimal solution of RLF. 

Noting that any instance of the RLF can be transformed into an instance of the well-known \textit{single-source capacitated facility location problem}, one can use the state-of-the-art exact algorithm given in \cite{Gadegaard2018} to solve RLF. However, we solve RLF by using a general purpose solver given that the time spent in solving RLF is negligible when compared to the total CPU time needed by our branch-and-cut algorithm. Finally, we only solve RLF when the support of the LP relaxation of the Benders reformulation changes from one iteration to the next at the root node, i.e., whenever $H(\bar{z}^{t-1}) \neq H(\bar{z}^t)$.


\subsubsection{Local Search Heuristic}\label{Subsub:Localsearch}

The local search heuristic is used to improve the initial solution obtained from the constructive heuristic. It consists of two phases. The first phase, denoted as the \textit{facility modification phase}, uses a variable neighborhood descent (VND) method to systematically explore five neighborhoods that modify both the set of hubs and assignment decisions. The second phase, denoted as the \textit{assignment modification phase}, solves a well-known MILP as an approximation to improve assignment decisions.

In what follows, solutions are represented by pairs of the form $s=(H, a)$ where $H \subseteq N$ denotes the set of opened facilities and $a: N \rightarrow H$ is the assignment mapping, i.e., $a(i)=k$ if customer $i\in N$ is assigned to facility $k \in H$. For any feasible assignment, $h_k$ denotes the available capacity of facility $k$, i.e., $h_k=b_k-\sum_{i:a(i)=k}d_i$.

During the facility modification phase, we use a VND method \citep{brimberg1995variable}. It is based on a systematic search in a set of $r$ neighborhoods, $\mathcal{N}_1,\mathcal{N}_2, \ldots, \mathcal{N}_r$. The VND works by performing a local search in a neighborhood $\mathcal{N}_1$ until a local optimal solution is found. After that, the algorithm switches to neighborhoods $\mathcal{N}_2, \ldots, \mathcal{N}_r$, sequentially, until an improved solution is found. Each time the search improves the best known solution, the procedure restarts using the neighborhood $\mathcal{N}_1$. Our implementation of the VND algorithm explores two types of neighborhood structures. The first type focuses on the reassignment of customers to open facilities, whereas the second type allows the set of open facilities to change. In all cases, we only consider movements that are feasible with respect to constraints \eqref{eq_Q_1}-\eqref{eq_Q_5}. 

The \textit{shift} neighborhood considers all solutions that can be obtained from the current one by changing the assignment of exactly one node, whereas the \textit{swap} neighborhood contains all solutions that differ from the current one in the assignment of two nodes. Let $s=(H, a)$ be the current solution, then
$$ \mathcal{N}_{1}(s)= \left\lbrace s'=(H, a'):\exists! i\in N, a'(i)\neq a(i) \right\rbrace ,$$
and
$$ \mathcal{N}_{2}(s)= \left\lbrace s'=(H, a'): \exists i_1, i_2,
a'(i_1)=a(i_2), a'(i_2)=a(i_1), a'(i)=a(i), \forall i \neq
i_1, i_2  \right\rbrace .$$

For exploring $\mathcal{N}_{1}$, we consider all pairs of the form $(i, j)$ where $a(j) \neq i $ and $h_i \geq d_j$. Also, for exploring $\mathcal{N}_{2}$ we consider all pairs of the form $(i_1, i_2)$ where $a(i_1) \neq a(i_2), h_{a(i_1)} + d_{i_1} \geq d_{i_2}$ and $h_{a(i_2)} + d_{i_2} \geq d_{i_1}$. In both cases we perform the first improving move.

We explore three additional neighborhood structures of the second type. They affect the current set of open facilities. The first one considers a subset of feasible solutions that are obtained from the current one by opening a new facility and by assigning some customers to it. That is, 
$$\mathcal{N}_{3}(s)\subset  \left\lbrace s'=(H', a'): H'=H\cup\left\{k\right\} ; \forall j, a'(j)=r \in H',  \sum_{j:a'(j)=r} d_j \leq b_r, \forall r \in H' \right\rbrace.$$

To explore $\mathcal{N}_{3}(s)$, all nodes $k \in N \setminus H$ are considered. Again, let $s=(H,a)$ denote the current solution, $a'(k)$ the new assignment and $\hat{h_r}$ the available capacity of hub $r$. Initially, $a'(p)=a(p)$ for all $p \in N$ and $\hat{h}_r=h_r$ for all $r \in H$. For each potential facility $k \in N \setminus H$, we consider nodes by decreasing order of their demand $d_i$. Node $j$ is reassigned to facility $k$ if $c_{jk} \leq c_{j a(j)}$ and $\hat{h}_k \geq d_j$. If node $j$ is reassigned to facility $k$ we update its assignment and the available capacity of facilities $k$ and $a(j)$.

The second neighborhood structure of the second type considers a subset of feasible solutions that are obtained from the current one by closing a facility and reassigning its assigned customers to other open facilities. That is,
$$\mathcal{N}_{4}(s)\subset \left\lbrace s'=(H', a'): H'=H \setminus k ; \forall j, a'(j)=r \in H', \sum_{j:a'(j)=r} d_j \leq b_r, \forall r \in H' \right\rbrace. $$

To explore $\mathcal{N}_{4}(s)$, all facilities $k \in H$ are considered. Once more, let $a'(k)$ denote the new assignment and $\hat{h_r}$ the available capacity of facility $r$. Initially, $a'(p)=a(p)$ for all $p \in N$ and $\hat{h}_r=h_r$, for all $r \in H$. For each open facility $k \in H$, we consider its assigned customers in decreasing order of their demand $d_i$. Customer $j$ is reassigned to facility $\hat{m}$, where $\hat{m}= argmin \left\{c_{jm} : \hat{h}_m - d_j \geq 0 , m \in H\setminus k \right\}$. If node $j$ is reassigned to facility $\hat{m}$, we update its assignment and the available capacity of facilities $\hat{m}$ and $a(j)$.

The last neighborhood of the second type considers a subset of feasible solutions that are obtained by simultaneously closing an open facility and opening a closed one.  That is, 
$$\mathcal{N}_5(s) \subset \left\lbrace s^{'} = (H^{'},a^{'} ):H^{'}=H \backslash \{m\}\cup\{i\},m \in H, \ i\in N \backslash H \right\rbrace. $$
To explore $\mathcal{N}_5(s)$, all nodes $i \in N \backslash H$ are considered, and a set of solutions is obtained from the current one by interchanging an open facility by a closed one while ensuring capacity is enough to meet demand, and reassigning all the customers to their closest open facility with available capacity.

Once the VND search terminates with a local optimal solution with respect to the above five neighborhoods, we proceed with the assignment modification phase. The objective of this second step is to intensify the search on the assignment decisions by solving a well-known combinatorial optimization problem that arises as a subproblem in the QC$p$LP. In particular, when the location of the facilities is given, the QC$p$LP reduces to a capacitated variant of the quadratic semi-assignment problem \citep{saito09,Rostami18}. Moreover, if we relax the quadratic term of \eqref{eq_Q_obj} for each set of open facilities $H \subseteq N$, we obtain the following \textit{generalized assignment problem} (GAP):

\begin{eqnarray}
(GAP) \quad \min\quad & & \sum_{i \in N}\sum_{k \in H} c_{ik}z_{ik} \label{eq_L_objGAP}\\
\mbox{s.t.} \quad & &    \sum_{k \in H}z_{ik} = 1 \qquad i \in N \label{eq_Q_1GAP}\\
& & \sum_{i\neq k} d_i z_{ik} \leq \bar{b}_{k}z_{kk}  \qquad  k \in H  \label{eq_Q_3GAP}\\
&  & z_{ik}\in \{0,1\}   \qquad i \in N, k \in H. \label{eq_Q_5GAP}
\end{eqnarray}

The solution to the GAP can be seen as an intensification mechanism in which we focus on improving the assignment decisions for all customers simultaneously, while taking explicitly into account the capacity restrictions at facilities but relaxing the quadratic term in the objective. In case the solution to the GAP for a given set $H$ obtained from the local optimal solution of the first step modifies some customers' assignments, we perform a simple VNS on the obtained solution in which only the shift and swap neighborhoods are explored. Similar to the MILP solved in the constructive phase, one can use the state-of-the-art exact algorithm given in \cite{avella10} to solve GAP. However, we solve GAP by using a general purpose solver given that the time spent in solving it is negligible compared to the total CPU time needed by our branch-and-cut algorithm.

\subsection{Model Reduction Techniques}\label{Subsec:ModelRed}
When solving large-scale optimization problems, reducing the size of the formulation without compromising optimality is crucial to reduce the computational time. Having a high quality feasible solution suggests the use of elimination tests and variable fixing procedures to significantly reduce the size of the formulation to be solved in the branch-and-cut algorithm. 

\subsubsection{Elimination Tests}\label{Subsub:Elim} 
First proposed in \citet{crowder1980}, elimination tests have been shown to be an effective tool for identifying variables that will not be in the optimal solution and can therefore be discarded. It uses linear programming principles to conclude that a chosen variable will have a value of zero at an optimal solution. The procedure is based on the following well-known proposition written in terms of our Benders reformulation $BR$.
\begin{proposition}
	Let $UB$ and $LB$ be an upper bound and an LP relaxation-based lower bound to $BR$, respectively,  and $rc_{k}$ be the reduced cost coefficient of $z_{kk}$, $k\in N$, at an optimal solution of any linear relaxation of $BR$ in which only a subset of Benders cuts is included.  If $LB + rc_{k} > UB$, then $z_{kk} = 0$ at an optimal solution to $\mbox{BR}$.
\end{proposition}

This proposition allows the elimination of a variable from the formulation without compromising on optimality. By verifying the condition for all non-basic location variables every time new constraints are added to the LP relaxation of BR, a significant number of variables may be eliminated. Note that we focus only on the location variables $z_{kk}$, $k\in H$ since the impact of removing it will lead to the elimination of all other related variables, $z_{ik}$, $\forall i \in N$, thereby reducing the size of BR by $|N|$.   

\subsubsection{Partial Enumeration}\label{Subsub:Partial}
Similar to the idea of strong branching \citep{Applegate95}, partial enumeration creates a set of what-if scenarios to measure the impact that fixing a variable $z_{kk}$, $k\in N$ to an arbitrary value $\pi$ has on the LP of the formulation. If the optimal solution value of the linear program with fixed variable $z_{kk}$ is larger than that of a known upper bound, then that variable will not have a value of $\pi$ at an optimal solution to BR. This method is particularly useful when the variable can only take two possible values as is the case with the binary location variables.

Since there are two possible what-if scenarios, one can create for each variable in the partial enumeration scheme, a simple rule to identify good candidates. If for a given location variable $z_{kk}$ $k\in N$ its value is less than or equal to 0.2, we then solve a linear programming problem with the added constraint $z_{kk}=1$. If the resulting optimal solution value is greater than the best upper bound known, we can then make $z_{kk}$ and its corresponding assignment variables equal to 0 (i.e., remove them from the formulation). We refer to this procedure as PE$_0$, which stands for partial enumeration to fix at 0.

If the value of $z_{kk}$ is greater than or equal to 0.8, we then solve a linear programming problem with the added constraint $z_{kk}=0$. If the resulting optimal solution value is greater than the best upper bound known, we then fix  $z_{kk}=1$. We refer to this procedure as PE$_1$, which stands for partial enumeration to fix at 1. 

Both the elimination test and partial enumeration procedure are performed as a part of a preprocessing phase in which the linear relaxation of BR is solved. This reduces the size of the mixed binary formulation used in the branch-and-cut algorithm, and of the modified transportation problems $\mbox{PPO$_{ij}$}$ solved to obtain Benders cuts. 



\subsection{Some Implementation Details of the Branch-and-Cut Algorithm for QC$p$LP}\label{Sec:CompleteAlgo}
Having explained each of the algorithmic enhancements of our procedure, we now present the complete exact algorithm which is divided into two parts. The first is the  root-processing routine, in which the linear programming relaxation of BR is solved and variables are dynamically eliminated and fixed. The second part is the branch-and-cut algorithm where the mixed binary program is solved to proven optimality. The transportation problem PPO$_{ij}$ is used to separate new Benders cuts at every node of the tree at a depth of multiples of $\gamma$.

The preprocessing phase is an iterative procedure in which the LP relaxation of BR without any Benders cuts is first solved. The mathheuristic is then called using the support (set of variables with non-zero solution value) of the LP relaxation to obtain a feasible solution and to potentially update the incumbent. We then perform the variable elimination procedure and, if the iteration number is divisible by $\phi$, we also perform the partial enumeration to close facilities (PE$_0$). We then proceed to generate a Benders cut by solving our modified transportation problem $\mbox{PPO$_{ij}$}$ for each $i,j\in N$ and add it to BR if the cut is violated by at least $\epsilon$. We then reoptimize the BR with its newly added constraint. If the resulting optimal solution value is larger than the previous by a margin of at least $\kappa$, the variable elimination procedure is executed again and the entire process is repeated. If the improvement is less than $\kappa$, we then execute the partial enumeration methods PE$_1$ followed by PE$_0$. The resulting formulation is then used in the mixed binary programming phase. 


Upon completing the preprocessing phase, the resulting mixed binary program BR is of a significantly smaller size than the original formulation. The reduced BR along with the previously generated Benders cuts are then inputted into a solver that executes a branch-and-cut process and infer additional general mixed integer cuts to strengthen the formulation. Up to $\Upsilon$ additional Benders cuts are obtained when exploring nodes at a depth of multiples of $\gamma$ of the branch-and-bound tree and are added only if the current fractional solution violates them by at least $\epsilon$. For correctness, the tolerance to add Benders cuts at integer solutions is set to 0. Finally, we also impose that in the branching framework, priority be given to location variables $z_{kk}$, to ensure the impact of opening and closing facilities is evaluated first. The source code can be downloaded from  \url{https://sites.google.com/view/carloszetina/Research/Publications}.

\section{Computational Experiments} \label{Sec:CompExp}

We next present the results of extensive computational experiments performed to \modif{assess} the performance of our branch-and-cut algorithm. Given that the star-star hub network design problem is a special case of hub location problems and that the instances are not publicly available, our experiments are based on the well-studied, well-benchmarked single-allocation hub location problems.  In the first part of the experiments, we focus on a comparison between our exact algorithm and other exact algorithms reported in the literature for four particular cases of QC$p$LP arising in hub location: the \textit{uncapacitated hub location problem with single assignments} (UHLPSA), the \textit{uncapacitated $p$-hub median problem with single assignments} (U$p$HMPSA), the \textit{capacitated hub location problem with single assignments} (CHLPSA), and the \textit{capacitated $p$-hub median problem with single assignments} (C$p$HMPSA). In the second part of the experiments, we test the robustness and limitations of our exact algorithm on large instances involving up to 1,000 nodes. All algorithms were coded in C using the callable library for CPLEX 12.9.0 and run on an Intel Xeon E5 2687W V3 processor at 3.10 GHz with 750 GB of RAM under Linux environment. The separation and addition of Benders optimality cuts within the branch-and-cut algorithm has been implemented via \textit{lazy callbacks} and \textit{user cut callbacks}. For a fair comparison, all use of CPLEX was limited to one thread and the traditional MIP search strategy.  The following parameter values were, in order of appearance used in our final implementation: $\phi=0.5$,  $\epsilon=100$, $\kappa=0.1$, $\Upsilon=2$, and $\gamma=10$.
 
\subsection{Comparison with Alternative Solution Algorithms}\label{Subsec:Compare}

We now present a comparison between our branch-and-cut algorithm for QC$p$LP and the most recently proposed exact algorithms for UHLPSA \citep{Meier2018,Ghaffarinasab2018}, U$p$HMPSA \citep{Ghaffarinasab2018}, and CHLPSA \citep{contreras2011b,Meier2018}. To the best of our knowledge, these are the state-of-the-art algorithms for solving the respecting problems to proven optimality. We also present the results of our algorithm for the C$p$HMPSA. However, to the best of our knowledge, there is no ad hoc exact algorithm in the literature for the C$p$HMPSA. All computational experiments in this section are performed using the well-known Australian Post (AP) set of instances. It consists of the Euclidean distances $d_{ij}$ between 200 postal districts in \modif{Sydney}, Australia, and of the values of $w_{ij}$ representing postal flows between pairs of postal districts. From this set of instances, we have selected those with $|N| \in \left\lbrace 100, 125, 150, 175, 200 \right\rbrace$, $p \in \left\lbrace 5, 10, 15, 20 \right\rbrace$ and with setup costs and capacities of the types loose (L) and tight (T) \citep[see][]{contreras2011b}. We can represent any instance of the above mentioned hub location problems as an instance of the QC$p$LP by setting the costs in the objective as $c_{ik} = \left(\chi O_i + \delta D_i \right)d_{ik}$ and $q_{ikjm} = \tau w_{ij}d_{km}$, where $O_i=\sum_{j \in N} w_{ij}$, $D_i=\sum_{j \in N} w_{ji}$, and $\chi$, $\tau$, $\delta$ represent the unit collection, transfer, and distribution costs, respectively. For the AP data set, these unit costs are set to $\chi=2$, $\tau=0.75$, and $\delta=3$. 

Our comparisons are done with \cite{Meier2018}, \cite{Ghaffarinasab2018}, and \cite{contreras2011b} which are, to the best of our knowledge, the state of the art ad hoc exact solvers for their respective problem variants.  \cite{Meier2018} use C$^{\#}$ to call Gurobi 6.0 running on a processor at 3.4 GHz while \cite{Ghaffarinasab2018} use Java to call CPLEX 12.6 on a processor at 3.30 GHz. Finally, \cite{contreras2011b} use C to develop a branch-and-price solver run on a processor at 2.33 GHz. 

Our first comparison is for the UHLPSA where there exist fixed set-up costs for locating facilities but both capacity and cardinality constraints are relaxed. We compare our algorithm to the results presented in \cite{Meier2018} and \cite{Ghaffarinasab2018}.  The detailed results of the comparison between the exact methods using the AP data set are provided in Table \ref{UHLPSAtable}. The first two columns are the instance name and optimal objective function value. These are then followed by the time in seconds reported by \cite{Ghaffarinasab2018} and \cite{Meier2018} and the information of our algorithm's performance. The column under the heading \textit{\%Dev heur} reports the percent deviation between the best solution obtained with the matheuristic presented in Section \ref{Subsec:Heuristic} and the optimal solution. The column \textit{\%fixed plants} gives the percent of facilities that were closed at the end of the partial enumeration phase. The column \textit{\%time root} provides the percent of time spent by the algorithm solving the root node, including the matheuristic, the elimination tests, and the partial enumeration phase. The last column reports the number of nodes explored in the enumeration tree.

\begin{table}[htbp]
	\TABLE
	{Comparison of state-of-the-art exact algorithms for the uncapacitated hub location problem with single assignments using the AP data set. \label{UHLPSAtable}}
	{\begin{tabular}{|rr|r|r|rrrrr|}
			\hline \up
			&       &       \multicolumn{1}{c|}{Ghaffarinasab}  & \multicolumn{1}{c|}{Meier and} & \multicolumn{5}{c|}{Branch-and-cut algorithm}  \\
			&       &        \multicolumn{1}{c|}{and Kara (2018)}  & \multicolumn{1}{c|}{Clausen (2018)} &  & \%Dev & \%fixed &  \%time & \multicolumn{1}{l|}{BB} \\
			Instance      & Opt.  & time(s) & time(s) & time(s) &  heur &  plants & root  & \multicolumn{1}{l|}{nodes} \\
			\hline \up
			100LT	&	238,016.28	&	6.91	&	82.42	&	23.76	&	0.00	&	97	&	68	&	70	\\
			100TT	&	305,097.95	&	3.24	&	60.23	&	2.32	&	1.52	&	98	&	68	&	0	\\
			125LT	&	227,949.00	&	43.30	&	411.46	&	42.89	&	0.00	&	98	&	37	&	122	\\
			125TT	&	258,839.68	&	16.08	&	188.05	&	8.54	&	0.00	&	98	&	21	&	9	\\
			150LT	&	225,450.09	&	107.24	&	1,259.22	&	76.05	&	0.00	&	99	&	39	&	172	\\
			150TT	&	234,778.74	&	26.30	&	478.38	&	59.06	&	0.00	&	99	&	5	&	227	\\
			175LT	&	227,655.38	&	188.17	&	2,044.77	&	89.86	&	0.00	&	99	&	70	&	51	\\
			175TT	&	247,876.80	&	44.55	&	1,639.21	&	69.84	&	0.00	&	99	&	8	&	180	\\
			200LT	&	233,802.98	&	68.18	&	5,493.47	&	290.87	&	0.00	&	99	&	65	&	276	\\
			200TT	&	272,188.11	&	1,399.53	&	20,292.35	&	102.27	&	0.46	&	97	&	29	&	136	\\
			\hline \up
			&	Geom. mean	&	45.15	&	782.85	&	41.96	&	0.00	&	98	&	31	&	0	\\
			&	Arith. mean	&	190.35	&	3,194.96	&	76.55	&	0.20	&	98	&	41	&	124	\\
			
			\hline 
	\end{tabular}}
	{}
\end{table}%

The results in Table \ref{UHLPSAtable} indicate that our exact algorithm \modif{is, on average, an order of magnitude faster than the cutting plane algorithm}  of \cite{Meier2018}, which explicitly \modif{requires the assumption} of Euclidean distances. \modif{On the other hand, our algorithm takes on average half of the computation time required by the iterative Benders decomposition algorithm of \cite{Ghaffarinasab2018}}. In addition, the matheuristic is capable of finding an optimal solution in 9 out of 10 instances. 

Our second comparison is for the U$p$HMPSA, where the fixed set-up costs for locating facilities are disregarded but a cardinality \modif{constraint} on the number of open hubs is imposed. Also, capacity limitations on the facilities are not considered. Once more, we compare our algorithm to the results presented in \cite{Meier2018} and \cite{Ghaffarinasab2018} for solving the U$p$HMPSA. The detailed results of the comparison between the exact methods using the AP data set are provided in Table \ref{UpHMPSAtable}.

   
\begin{table}[ht!]
	\TABLE
	{Comparison of state-of-the-art exact algorithms for the uncapacitated p-hub median problem with single assignments using the AP data set. \label{UpHMPSAtable}}
	{\begin{tabular}{|rrr|r|r|rrrrr|}
			\hline \up
	     	&       &       & \multicolumn{1}{c|}{Ghaffarinasab}  & \multicolumn{1}{c|}{Meier and} & \multicolumn{5}{c|}{Branch-and-cut algorithm}  \\
			&       &       & \multicolumn{1}{c|}{and Kara (2018)}  & \multicolumn{1}{c|}{Clausen (2018)} &  & \%Dev & \%fixed &  \%time & \multicolumn{1}{l|}{BB} \\
   Instance & p     & Opt.  & time(s) & time(s) & time(s) &  heur &  plants & root  & \multicolumn{1}{l|}{nodes} \\
\hline \up
		100	&	5	&	136,929.44	&	313.80	&	356.39	&	23.26	&	0.00	&	86	&	82	&	32	\\
		100	&	10	&	106,469.57	&	109.18	&	32.22	&	34.33	&	0.00	&	82	&	84	&	14	\\
		100	&	15	&	90,533.52	&	144.40	&	85.95	&	41.92	&	0.00	&	67	&	66	&	79	\\
		100	&	20	&	80,270.96	&	61.47	&	33.54	&	26.35	&	0.00	&	72	&	82	&	10	\\
		125	&	5	&	137,175.68	&	1,286.48	&	1,104.31	&	79.83	&	0.04	&	88	&	58	&	111	\\
		125	&	10	&	107,092.09	&	414.08	&	184.62	&	69.69	&	0.00	&	81	&	89	&	13	\\
		125	&	15	&	91,494.56	&	1,271.86	&	465.08	&	81.34	&	0.00	&	69	&	65	&	132	\\
		125	&	20	&	81,471.65	&	213.76	&	111.03	&	107.45	&	0.00	&	67	&	63	&	149	\\
		150	&	5	&	137,425.90	&	2,989.83	&	1,474.76	&	153.27	&	0.00	&	91	&	72	&	178	\\
		150	&	10	&	107,478.12	&	1,148.01	&	412.53	&	159.79	&	0.00	&	92	&	82	&	68	\\
		150	&	15	&	92,050.58	&	1,695.15	&	405.57	&	256.22	&	0.03	&	68	&	52	&	142	\\
		150	&	20	&	82,229.39	&	531.99	&	185.34	&	195.07	&	0.00	&	78	&	70	&	114	\\
		175	&	5	&	139,354.51	&	31,347.15	&	10,699.58	&	535.70	&	0.00	&	92	&	54	&	342	\\
		175	&	10	&	109,744.35	&	10,551.64	&	3,023.02	&	895.67	&	0.61	&	57	&	41	&	236	\\
		175	&	15	&	94,123.66	&	19,602.93	&	8,143.73	&	627.45	&	0.04	&	66	&	42	&	402	\\
		175	&	20	&	83,843.59	&	1,778.11	&	271.03	&	429.19	&	0.00	&	79	&	61	&	167	\\
		200	&	5	&	140,062.65	&	127,546.79	&	17,628.38	&	1,156.65	&	0.07	&	79	&	50	&	367	\\
		200	&	10	&	110,147.66	&	46,706.90	&	4,957.31	&	1,302.66	&	0.00	&	83	&	39	&	1036	\\
		200	&	15	&	94,459.20	&	26,640.56	&	1,107.48	&	769.89	&	0.21	&	66	&	52	&	151	\\
		200	&	20	&	84,955.36	&	27,224.48	&	526.08	&	1,141.73	&	0.00	&	75	&	37	&	1007	\\
		\hline \up 
	&		&	Geom. mean	&	2,195.82	&	613.97	&	198.46	&	0.00	&	76	&	60	&	122.88	\\
	&		&	Arith. mean	&	15,078.93	&	2,560.40	&	404.37	&	0.05	&	77	&	62	&	237.5	\\
		\hline 
	\end{tabular}}
	{}
\end{table}%

The results in Table \ref{UpHMPSAtable} show once more that our exact algorithm \modif{is on average three and two orders of magnitude faster than those of \cite{Ghaffarinasab2018} and \cite{Meier2018}, respectively}. The matheuristic is capable of finding an optimal solution in 14 out of 20 instances, and the average deviation is only 0.05\%. It is worth noting that the U$p$HMPSA seems to be more difficult to solve as compared to the UHLPSA. The largest CPU time for the UHLPSA was about six minutes (200 nodes with loose setup costs) whereas for the U$p$HMPSA was one hour (200 nodes and $p=20$).

Our third comparison is for the CHLPSA, where the fixed set-up costs and capacities of facilities are taken into account but the cardinality constraint is disregarded. We compare our algorithm to the results presented in \cite{Meier2018} using a branch-and-cut algorithm and \cite{contreras2011b} using a branch-and-price algorithm for solving the CHLPSA. The detailed results of the comparison between the exact methods using the AP data set are provided in Table \ref{CHLPSAtable}. 

\begin{table}[ht!]
	\TABLE
	{Comparison of state-of-the-art exact algorithms for the capacitated hub location problem with single assignments using the AP data set. \label{CHLPSAtable}}
	{\begin{tabular}{|rr|r|r|rrrrr|}
			\hline \up
			&       &       \multicolumn{1}{c|}{Meier and}  & \multicolumn{1}{c|}{Contreras} & \multicolumn{5}{c|}{Branch-and-cut algorithm}  \\
			&       &        \multicolumn{1}{c|}{Clausen (2018)}  & \multicolumn{1}{c|}{ et al. (2011d)} &  & \%Dev & \%fixed &  \%time & \multicolumn{1}{l|}{BB} \\
			Instance      & Opt.  & time(s) & time(s) & time(s) &  heur &  plants & root  & \multicolumn{1}{l|}{nodes} \\
			\hline \up
			100LL	&	246,713.97	&	176.31	&	459.89	&	19.11	&	0.00	&	93	&	54	&	68	\\
			100LT	&	256,155.33	&	141.23	&	347.95	&	27.88	&	0.06	&	95	&	29	&	155	\\
			100TL	&	362,950.09	&	174.03	&	124.92	&	40.49	&	0.59	&	92	&	14	&	242	\\
			100TT	&	474,068.96	&	1,122.14	&	328.11	&	95.35	&	0.37	&	93	&	6	&	943	\\
			125LL	&	239,889.33	&	415.41	&	1,650.57	&	45.15	&	0.02	&	93	&	52	&	101	\\
			125LT	&	251,259.16	&	279.53	&	552.99	&	45.11	&	0.02	&	97	&	41	&	299	\\
			125TL	&	246,486.69	&	53.24	&	41.22	&	59.04	&	0.00	&	98	&	4	&	302	\\
			125TT	&	291,807.35	&	109.19	&	322.73	&	27.52	&	0.00	&	97	&	12	&	194	\\
			150LL	&	234,765.44	&	153.85	&	3,347.21	&	53.01	&	0.35	&	97	&	82	&	8	\\
			150LT	&	249,797.49	&	time out	&	11,818.19	&	84.61	&	0.00	&	92	&	50	&	293	\\
			150TL	&	262,543.08	&	1,419.80	&	1,114.95	&	73.29	&	0.31	&	98	&	13	&	393	\\
			150TT	&	322,976.47	&	3,243.56	&	4,299.28	&	192.74	&	0.00	&	92	&	7	&	1124	\\
			175LL	&	227,997.58	&	1,010.55	&	3,418.10	&	119.05	&	0.00	&	98	&	56	&	150	\\
			175LT	&	251,540.80	&	2,196.45	&	12,408.05	&	131.39	&	0.00	&	95	&	61	&	138	\\
			175TL	&	244,860.41	&	152.07	&	256.60	&	54.48	&	0.00	&	98	&	11	&	138	\\
			175TT	&	312,193.78	&	17,634.54	&	4,886.88	&	223.70	&	0.07	&	94	&	7	&	764	\\
			200LL	&	231,069.50	&	725.40	&	5,813.00	&	210.83	&	0.00	&	98	&	54	&	241	\\
			200LT	&	267,218.35	&	7,460.18	&	45,874.73	&	980.45	&	0.71	&	86	&	19	&	1441	\\
			200TL	&	273,443.81	&	902.96	&	869.67	&	163.37	&	0.00	&	98	&	8	&	412	\\
			200TT	&	290,582.04	&	4,117.31	&	3,211.04	&	346.75	&	0.11	&	98	&	7	&	804	\\
			\hline \up
			&	Geom. mean	&	470.16	&	1,376.70	&	88.96	&	0.00	&	95	&	20	&	252	\\
			&	Arith. mean	&	2,183.57	&	5,057.30	&	149.67	&	0.13	&	95	&	29	&	411	\\
			
			\hline 
	\end{tabular}}
	{}
\end{table}%

In the case of the CHLPSA, our algorithm \modif{is significantly faster than} the others. In some instances it is up to three orders-of-magnitude faster, e.g., 150LT which was not solved in 12 hours by  \cite{Meier2018}. The matheuristic provides an optimal solution for 13 out of 20 instances, while the average deviation is only 0.13\%. The most time consuming instance (200LT) was solved in about eight minutes.

Our last series of experiments in this section is for the C$p$HMPSA, where the fixed set-up costs are disregarded but capacities and cardinality constraints of facilities are taken into account. To the best of our knowledge, there is \modif{no} ad hoc algorithm for solving this hub location variant in the literature. Therefore, in Table \ref{CpHMPSAtable}, we only report the results obtained with our algorithm using the AP data set. 

\begin{table}[ht!]
	\TABLE
	{Results of branch-and-cut algorithm for the capacitated $p$-hub median problem with single assignments using the AP data set. \label{CpHMPSAtable}}
	{\begin{tabular}{|rrr|rrrrr|}
			\hline \up
			&    &   &   \multicolumn{5}{c|}{Branch-and-cut algorithm}  \\
			&    &   &   & \%Dev & \%fixed &  \%time & \multicolumn{1}{l|}{BB} \\
			Instance   & $p$   & Opt.  & time(s) & heur &  plants & root  & \multicolumn{1}{l|}{nodes} \\
			\hline \up
		100	&	5	&	137,232.52	&	26.08	&	0.04	&	88	&	70	&	30	\\
		100	&	10	&	107,207.72	&	34.69	&	0.00	&	80	&	81	&	12	\\
		100	&	15	&	91,283.34	&	40.10	&	0.00	&	65	&	61	&	63	\\
		100	&	20	&	81,034.59	&	25.73	&	0.00	&	71	&	76	&	22	\\
		125	&	5	&	137,175.68	&	118.32	&	0.00	&	86	&	39	&	313	\\
		125	&	10	&	107,092.09	&	65.62	&	0.00	&	90	&	90	&	18	\\
		125	&	15	&	91,494.56	&	84.85	&	0.00	&	69	&	62	&	96	\\
		125	&	20	&	81,471.65	&	127.93	&	0.00	&	67	&	53	&	188	\\
		150	&	5	&	137,425.90	&	152.15	&	0.00	&	90	&	74	&	117	\\
		150	&	10	&	107,478.12	&	159.42	&	0.00	&	92	&	83	&	82	\\
		150	&	15	&	92,050.58	&	278.94	&	0.03	&	67	&	51	&	163	\\
		150	&	20	&	82,229.39	&	203.05	&	0.00	&	75	&	69	&	82	\\
		175	&	5	&	139,354.51	&	544.00	&	0.00	&	92	&	55	&	462	\\
		175	&	10	&	109,744.35	&	856.88	&	0.61	&	57	&	45	&	205	\\
		175	&	15	&	94,123.66	&	652.80	&	0.01	&	67	&	42	&	343	\\
		175	&	20	&	83,843.59	&	474.33	&	0.00	&	79	&	61	&	182	\\
		200	&	5	&	140,062.65	&	1,018.33	&	0.08	&	82	&	62	&	260	\\
		200	&	10	&	110,147.66	&	1,274.57	&	0.29	&	70	&	44	&	616	\\
		200	&	15	&	94,459.20	&	911.05	&	0.21	&	66	&	46	&	176	\\
		200	&	20	&	84,955.37	&	1,720.54	&	0.00	&	75	&	28	&	1084	\\
		\hline \up
		&		&	Geom. mean	&	211.39	&	0.00	&	75.61	&	57.32	&	126.91	\\
		&		&	Arithm. mean	&	438.47	&	0.06	&	76.32	&	59.61	&	225.7	\\

			\hline 
	\end{tabular}}
	{}
\end{table}%

The results of Table \ref{CpHMPSAtable} show that our exact algorithm can optimally solve all considered instances in less than one hour. In fact, the most time consuming instance (200 nodes and $p=20$) took only 30 minutes to be solved. The matheuristic is capable of finding an optimal solution in 13 out of 20 instances, and the average deviation is only 0.06\%. \modif{ It is worth noting that the C$p$HMPSA seems to be more difficult to solve as compared to the CHLPSA. The largest CPU time for the CHLPSA is about eight minutes (200 nodes with loose setup costs and tight capacities) whereas for the C$p$HMPSA is 30 minutes (200 nodes and $p=20$)}. 

\subsection{Results for Larger Instances}\label{Subsec:Largeinstances}
As seen in Section \ref{Subsec:Compare}, our algorithm scales well to larger instances while either being competitive or significantly faster than \modif{tailored exact algorithms that exploit problem-specific structure for all considered problem variants}. This leads to the question of what are the size limits of our proposed algorithm. In this section, we present the results of solving the large-scale instances first presented in \cite{contrerasb} as a single-assignment variant. During this last set of experiments, we focus on the two variants in which our algorithm performed the best for the AP data set: the UHLPSA and CHLPSA.

\cite{contrerasb} introduce a new data set containing three different sets of instances with diverse structural characteristics in the flow network. They consider different levels of magnitude for the amount of flow originating at a given node to obtain three different sets of nodes: low-level (LL) nodes, medium-level (ML) nodes, and high-level (HL) nodes. The total outgoing flow of LL, ML, and HL nodes lies in the interval $[1,10]$, $[10,100]$, and $[100,1,000]$, respectively. In this section, we use the first set of instances, called \textit{Set I}, in which the number of HL, ML, and LL nodes is 2\%, 38\%, and 60\% of the total number of nodes, respectively. We use the instance generation code of \cite{contrerasb} to generate instances with $N=$ 250, 300, 350, 400, 450, 500, 550, 600, 650, 700, 750, 800, 850, 900, 950, and 1,000. To generate setup costs and capacities, we use the same procedure as described in \cite{contrerasb,contreras11b}.



Tables \ref{LargeUncap} and \ref{LargCap} present the detailed results of our branch-and-cut algorithm for the UHLPSA and CHLPSA, respectively, using the \textit{Set I} instances. In the case of the UHLPSA, instances with up to 700 nodes can be optimally solved in less than half of a day of CPU time. The largest-size instances with 750 to 1,000 were optimally solved in CPU times of less than three days. From our observations, these results are remarkable given that the number of variables in the largest instance with 1,000 nodes require 499,501,000,000 variables in the $RL_2$ formulation. Similar results are obtained for the more challenging CHLPSA, in which instances with up to 800 nodes can be solved in less than one day of CPU time. The largest instances from 850 to 1,000 nodes require two to three days to prove optimality. Moreover, the matheuristic is capable of finding the optimal solution to 13 out of 15 uncapacitated instances and in 14 out of 20 for the capacitated case, with an average deviation of only 0.01\%.  

\begin{table}[ht!]
	\TABLE
	{Results for the uncapacitated hub location problem with single assignments using \textit{Set I} data set. \label{LargeUncap}}
	{\begin{tabular}{|rr|rrrrr|}
			\hline \up
			&      &   \multicolumn{5}{c|}{Branch-and-cut algorithm}  \\
			&      &   & \%Dev & \%fixed &  \%time & \multicolumn{1}{l|}{BB} \\
			Instance     & Opt.  & time(s) & heur &  plants & root  & \multicolumn{1}{l|}{nodes} \\
			\hline \up
			250	&	7,349,579.02	&	451.85	&	0.03	&	94	&	42	&	340	\\
			300	&	7,499,305.16	&	1,103.00	&	0	&	96	&	49	&	687	\\
			350	&	9,361,597.79	&	2,203.65	&	0	&	98	&	59	&	745	\\
			400	&	22,060,422.84	&	2,821.72	&	0	&	95	&	66	&	699	\\
			450	&	28,765,424.39	&	6,884.78	&	0	&	96	&	62	&	1560	\\
			500	&	35,679,997.04	&	13,536.08	&	0.01	&	95	&	31	&	4007	\\
			600	&	46,396,548.35	&	18,716.31	&	0	&	99	&	83	&	1489	\\
			650	&	50,183,822.62	&	31,409.75	&	0	&	99	&	84	&	1679	\\
			700	&	59,187,662.18	&	43,505.48	&	0	&	98	&	75	&	3767	\\
			750	&	67,824,339.65	&	63,430.23	&	0	&	98	&	83	&	3293	\\
			800	&	74,897,428.36	&	69,898.70	&	0	&	98	&	88	&	1493	\\
			850	&	84,429,662.33	&	142,740.82	&	0	&	97	&	66	&	6296	\\
			900	&	91,522,770.07	&	99,430.21	&	0	&	98	&	83	&	2537	\\
			950	&	100,422,604.52	&	182,577.62	&	0	&	99	&	77	&	7277	\\
			1000	&	111,339,940.06	&	226,876.92	&	0	&	98	&	62	&	7248	\\
			\hline \up
			&	Geom. mean	&	19,303.40	&	0.00	&	97	&	65	&	1987	\\
			&	Arithm. mean	&	60,372.47	&	0.00	&	97	&	67	&	2874	\\
			\hline 
	\end{tabular}}
	{}
\end{table}%

\begin{table}[ht!]
	\TABLE
	{Results for the large-scale capacitated hub location problem with single assignments using \textit{Set I} data set. \label{LargCap}}
	{\begin{tabular}{|rr|rrrrr|}
			\hline \up
			&      &   \multicolumn{5}{c|}{Branch-and-cut algorithm}  \\
			&      &   & \%Dev & \%fixed &  \%time & \multicolumn{1}{l|}{BB} \\
			Instance     & Opt.  & time(s) & heur &  plants & root  & \multicolumn{1}{l|}{nodes} \\
			\hline \up
		250	&	7,880,912.40	&	439.14	&	0.04	&	94	&	58	&	207	\\
		300	&	7,611,876.55	&	1,194.01	&	0.00	&	96	&	43	&	991	\\
		350	&	9,445,459.03	&	6,524.52	&	0.00	&	95	&	18	&	5338	\\
		400	&	22,828,621.74	&	5,976.75	&	0.05	&	94	&	29	&	2259	\\
		450	&	29,715,150.35	&	20,463.62	&	0.20	&	92	&	19	&	4968	\\
		500	&	36,500,689.64	&	5,221.38	&	0.00	&	98	&	81	&	498	\\
		600	&	47,878,571.44	&	30,162.95	&	0.02	&	96	&	56	&	2385	\\
		650	&	51,115,098.98	&	57,811.47	&	0.00	&	98	&	48	&	10082	\\
		700	&	60,386,958.62	&	50,582.63	&	0.00	&	97	&	54	&	5968	\\
		750	&	68,382,804.86	&	49,143.18	&	0.01	&	99	&	78	&	3002	\\
		800	&	75,011,490.17	&	60,505.13	&	0.00	&	99	&	86	&	2065	\\
		850	&	85,547,053.64	&	102,411.16	&	0.01	&	99	&	80	&	4037	\\
		900	&	92,987,214.82	&	111,830.99	&	0.00	&	99	&	75	&	5260	\\
		950	&	101,245,457.94	&	212,247.85	&	0.00	&	99	&	59	&	14463	\\
		1000	&	113,144,817.97	&	227,549.00	&	0.00	&	99	&	62	&	8719	\\
		\hline \up
		&	Geom. mean	&	23,266.50	&	0.00	&	97	&	51	&	3003	\\
		&	Arithm. mean	&	62,804.25	&	0.02	&	97	&	56	&	4683	\\
			\hline 
	\end{tabular}}
	{}
\end{table}%

\section{Conclusions} \label{Sec:Conclusions}
In this paper, we have studied a general class of non-convex quadratic capacitated $p$-location problems with single assignments. The quadratic term in the objective function accounts for the interaction cost between facilities and customer assignments. A binary quadratic program was first linearized by applying a reformulation-linearization technique and the associated additional variables were then projected out using Benders decomposition. We proposed an exact branch-and-cut algorithm that incorporated several algorithmic refinements such as stabilized Benders Pareto-optimal cuts, a matheuristic, and variable elimination and partial enumeration procedures. Extensive computational experiments on a set of existing and new large-scale instances with up to 1,000 nodes have clearly confirmed the efficiency and robustness of our algorithm in solving four different particular cases of the studied problem. To the best of our knowledge, the new instances are by far the largest and most difficult ever solved for any type of hub location problem.

\ACKNOWLEDGMENT{%
This research was partially funded by the the Canadian Natural Sciences and Engineering research Council under grants 2018-06704, 2017-01746, and by a Research \& Publication Grant, Indian Institute of Management Ahmedabad.  The second author would also like to acknowledge the support and productive research environment provided by Andrea Lodi at the Canada Excellence Research Chair in Data Science for Real-Time Decision-Making as well as IVADO Labs. This support is gratefully acknowledged.
}

\bibliographystyle{myinforms2014} 
\bibliography{BibliographyThesis} 


\end{document}